\documentclass{birkart}
\usepackage{amsmath}
\usepackage{graphicx}
\usepackage{color}
\input amssym.def
\input amssym

 \newtheorem{thm}{Theorem}[section]

 \newtheorem{cor}[thm]{Corollary}

 \newtheorem{lem}[thm]{Lemma}

 \newtheorem{prop}[thm]{Proposition}

 \theoremstyle{definition}

 \theoremstyle{remark}

 \numberwithin{equation}{section}

\DeclareSymbolFont{AMSb}{U}{msb}{m}{n}
\DeclareMathSymbol{\N}{\mathbin}{AMSb}{"4E}
\DeclareMathSymbol{\Z}{\mathbin}{AMSb}{"5A}
\DeclareMathSymbol{\R}{\mathbin}{AMSb}{"52}
\DeclareMathSymbol{\Q}{\mathbin}{AMSb}{"51}
\DeclareMathSymbol{\I}{\mathbin}{AMSb}{"49}
\DeclareMathSymbol{\C}{\mathbin}{AMSb}{"43}

\def\M{\operatorname{M}}

\def\({\left(}
\def\){\right)}
\def\[{\left[}
\def\]{\right]}

\def\l{\lambda}

\def\o{\omega}

\def\g{\gamma}
\def\d{\delta}
\def\e{\epsilon}
\def\p{\partial}

\def\s{\sigma}

\def\G{\Gamma}

\def\i{\infty}
\def\O{\Omega}
\begin{document}

\title[Plateau's Problem]
  {On Plateau's Problem for Soap Films with a Bound on Energy}

\author[J. Harrison]{Jenny Harrison}

\address{  Department of Mathematics\\ University of California, Berkeley\\Berkeley, CA\\
94705\\USA}

\email{harrison@math.berkeley.edu}

%\subjclass{Primary 28A75; Secondary 49Q15}


 \date{July 23,   2003}

\begin{abstract}  We prove existence and a.e. regularity of an area minimizing soap film with a bound on energy spanning a given Jordan curve in $R^3$.

\end{abstract}

\maketitle

 \section{Introduction}   Given a simple, closed curve in
three-space is there a surface with minimal area spanning it?  Solutions to the problem of Plateau depend, of course, on the class of spanning surfaces permitted.    
Douglas \cite{Douglas} won the first Fields' medal for his proof of the existence of an area minimizing mapping of the 2-disk whose image spans the Jordan curve.  Regularity took
many years to establish \cite{Osserman} and some aspects are still unresolved.  Federer and Fleming's solutions  \cite{FF}  are area minimizing in the class of integral currents.   
Two years later Fleming \cite{Fleming1} proved any such solution  is an embedded, orientable surface, smooth away from its boundary.   
None of these solutions consider soap films that arise in nature such as Moebius strips or films with triple branching.   Almgren \cite{A} invented varifolds to treat
soap films, but the lack of a natural boundary operator slowed progress and regularity was never proved. 

Plateau observed that soap
films have only two possible kinds of branching:  (1) three sheets of surface meeting at
$120^{\circ}$ angles along a curve and (2) four such curves meeting at approximately  $109^{\circ}$ angles at a point \cite{Plateau}.   In \cite{soapfilm} the author provides
models for surfaces called {\em flat dipolyhedra} that  model all such films, orientable or nonorientable, as well as the surfaces considered in \cite{Douglas} and \cite{FF}.    Flat dipolyhedra take
advantage of the fact that soap films are actually two films essentially occupying the same space, but are not cancelling.     There is a natural boundary operator of
flat $k$-dimensional dipolyhedra into
flat $(k-1)$-dimensional dipolyhedra which relies on cohesion of the soap film structure supported by a geometric version of Cartan's magic formula.

The {\em energy} of a flat dipolyhedron is defined to be the length  of the singular branched
set of
$A$ plus the surface area of $A$.  In this paper we prove the existence of a  surface spanning a Lipschitz Jordan curve that is area minimizing in the space of flat dipolyhedra $A$ with  energy bounded by a  fixed
constant.   

Complete soap film regularity remains an open question, although we  prove here that the solution  is a smooth surface away from its branched set which is a   union of Lipschitz Jordan
curves of finite total length.

 \section{Preliminaries}
\subsection*{Chainlets} A {\em $k$-cell}  $\s$ in $\R^k$ is defined to be a finite intersection of $k$-dimensional half spaces.  A {\em $k$-cell in $\R^n$} is
a
$k$-cell in a $k$-dimensional subspace of $\R^n.$   Each $k$-cell is assumed to be oriented.  The  support  of a $k$-cell $\s$ is denoted by $|\s|.$  An {\em integral cellular
$k$-chain} is a formal sum
$\sum a_i
\s_i$  where
$a_i \in \Z $ and the $\s_i$ are $k$-cells.  Equate two cellular chains
$S_1 \sim S_2$   if and only if they have a common cellular subdivision which is nonoverlapping.  A polyhedral $k$-chain is an equivalence class of cellular
$k$-chains. The {\em mass} of a $k$-cell $\s$ is its $k$-dimensional Hausdorff measure $M(\s)
= \mathcal{H}_k(|\s|).$  If
$P$ is a polyhedral $k$-chain represented by a nonoverlapping $k$-chain $
\sum a_i \s_i $,   its    mass is defined by $M(P) = \sum |a_i|M(\s_i).$  

Let $G$ be an abelian group with a translation invariant metric making it a complete metric space.      Let $|g|$ denote the distance between $g \in
G$ and the group identity.       If $H$ is a closed subgroup of $G$, we use the quotient metric $|\bar{g}| = \inf\{|g|: g \in
\bar{g}\}.$  If $G = Z$ then $|g|$ denotes the absolute value.  The group $Z_p = \Z /p\Z$ is of special interest and we give it the quotient metric.   Let ${\mathbf P}_k(G) =
G
\otimes {\mathbf P}_k(\Z).$  This is the group of polyhedral $k$-chains with coefficients in $G$.   If $P \in {\mathbf P}_k(G) $ then $P = \sum g_i \s_i$ where the
$\s_i$ are nonoverlapping $k$-cells.  If $P \in {\mathbf P}_k(G)$, define $M(P) = \sum |g_i| M(\s_i).$  Whitney's flat norm on polyhedra is defined by  $$M_{\flat}(P) = \inf\{M(Q) +
M(R): D = Q + \p R, Q
\in {\mathbf P}_k(G), R \in {\mathbf P}_{k+1}(G)\}.$$ The completion of the group $ {\mathbf P}_k(G) $ with the norm $M_{\flat}$ is denoted $\mathcal{P}_k(G)$.  Its elements are called
{\em flat chains} with coefficients in $G$. The {\em support} of a flat chain $A$ is well defined (see \cite{whitney}) and is denoted $|A|.$  Define subgroups

$$\begin{aligned} \mathcal{M}_k(G) &= \{A \in  {\mathcal P}_k(G): M(A) < \i\}\\ \mathcal{N}_k(G) &= \{A \in  {\mathcal P}_k(G):
M(A)+ M(\p A) < \i\} \mbox{ and } \\\mathcal{N}_k^0(G) &= \{A \in  {\mathcal N}_k(G):
|A| \mbox{ is compact}\}.\end{aligned}$$

If
$\sigma$ is a cell and $v$ a vector in $\R^n$ then $T_{v}\sigma$ denotes the translation of
$\sigma$ through  
$v$.  A {\em
$1$-multicell} is a cellular chain of the form
$\sigma^1 = \sigma^0-T_{v_1} \sigma^0 $  where $\sigma^0$ is a cell and $v_1$ is a vector.
 Given a    vector    $v_j$ and a
$2^{j-1}$-multicell $\sigma^{j-1}$, define the {\em   $2^j$-multicell}
$\sigma^j$ as the cellular chain $\sigma^j = \sigma^{j-1} -T_{v_j}\sigma^{j-1}.$  Thus
$\sigma^j$ is generated by  vectors $v_1, \dots ,v_j$ and a cell $\sigma^0.$
  An {\em  integral $2^j$-multicellular  chain in $\R^n$} is a formal sum of
$2^j$-multicells, $S^j =
\sum_{i=1}^n a_i \sigma_i^j $ with coefficients $a_i \in \Z$.  
 Let $j \ge 1.$  Given a $2^j$-multicell $\sigma^j$ generated by a cell $\sigma^0$ and vectors   $v_1, \cdots, 
v_j$,  define  
$\|\sigma^j\|_j = M(\sigma^0)|v_1| |v_2 |\cdots  |v_j| $  where $|v|$ denotes the norm of a vector $v \in \R^n$.  For consistency of  notation define  $\|\sigma^0\|_0 = M(\sigma^0). $ 
 For $S^j = \sum a_i \sigma_i^j$ define $\|S^j\|_{j} = \sum_{i=1}^n
|a_i|\|\sigma_i^j\|_{j}.$

Let $ {\mathbf P}_k(\Z)$ denote polyhedral $k$-chains in $\R^n$.  Suppose $P  \in {\mathbf P}_k $   and $r \in
\Z^+$.  For $r = 0$ define
$|P|^{\natural_0} = M(P).$
For $r \ge 1$ define the {\em $r$-natural} norm $$ |P|^{\natural_r} = \inf\left\{ \sum_{s=0}^r \|S^j\|_j   +
  |C|^{\natural_{r-1}}\right\}
$$
where the infimum is taken over all  decompositions $P
= \sum_{s=0}^{r} [S^j]  + \partial C $ with $S^j$ a $2^j$ multicellular $k$-chain and $C$ a
polyhedral
$(k+1)$-chain.    The group of polyhedral $k$-chains ${\mathbf P}_k(\Z)$ completed with
the norm $|\quad|^{\natural_r}$ is denoted ${\mathcal
N}_{k}^{\natural_r}(\Z).$ Elements of ${\mathcal N}_{k}^{\natural_r}(\Z)$ are
called {\em $k$-dimensional chainlets of class $N^r$.}   The boundary operator $$\p:{\mathcal
N}_{k}^{\natural_r}(\Z) \to {\mathcal
N}_{k-1}^{\natural_{r+1}}(\Z)$$ is naturally defined and continuous.

Let $\s$ be a $k$-cell supported in $\R^n$ and set $v =  (0,0,\cdots, 0,1) \in \R^{n+1}.$    In \cite{soapfilm} {\em mass cells} are defined as limits in the $1$-natural norm $\mu
\s =
\lim_{h \to 0} \frac{\s \times hv}{|h|}.$  Integral mass chains are finite sums of mass cells with integer coefficients.
   Although the support of a mass cell is $k$-dimensional, as a current it is $(k+1)$-dimensional.   {\em Dipole cells} are defined by $\d \s = \p
\mu
\s +
\mu
\p \s.$  Dipole chains $S$ are  finite sums of dipole cells.  The {\em  support} of a mass chain  $T = \sum a_i \mu \s_i$ is  defined by $|T |  =
\cup |\s_i|$ and the {\em support} of a dipole chain $S = \sum  a_i \d \s_i$ is defined to be $|S| = \cup |\s_i|.$    Define the {\em mass} of a mass chain $T$   by 
$M(T)   = \sum |a_i|M(\s_i)$ and the {\em weight} of a  dipole chain $S$ by $W(S) = \sum |a_i|M(\s_i)$  where the $\{\s_i\}$ are nonoverlapping.   Spaces of mass and
dipole chains are denoted
${\mathbf T}_k(\Z)$ and ${\mathbf S}_k(\Z)$, resp.  Let $G$ be an abelian group.  By taking the tensor product with $G$, we may define   spaces of mass
and dipole chains with coefficients in $G$ and denote them by  
${\mathbf T}_k(G)$ and ${\mathbf S}_k(G)$, resp.
    Define  the space of {\em dipolyhedra} as the direct sum 
 ${\mathbf D}_k(G) = {\mathbf S}_k(G) \oplus {\mathbf T}_k(G).$
 If $D = S_D + T_D = S+T \in {\mathbf D}_k(G)$ define  the {\em energy} of $D$ by
$E(D) = W(S) + M(T).$  Finally define the $E_{\flat}$ norm on the space of $k$-dimensional dipolyhedra  by 
$$E_{\flat}(D) = \inf\{E(Q) + E(R): D = Q + \p R, Q \in {\mathbf D}_k(G), R \in {\mathbf D}_{k+1}(G)\}.$$ The completion of the space of dipolyhedra ${\mathbf
D}_k(G)$ with the
$E_{\flat}$ norm is an abelian group denoted ${\mathcal D}_k(G).$  The boundary operator is continuous in the $E_{\flat}$ norm and satisfies $ E_{\flat}(\p A) \le
E_{\flat}(A).$   An element  $A \in {\mathcal D}_k(G)$ is called a {\em flat dipolyhedron} with coefficients in $G$.   In \cite{soapfilm} it is shown that    
weight, mass and energy  are well defined and  lower semi-continuous  in
${\mathcal D}_k(G)$.  
The operators $\d$ and $\mu$ satisfy
$E_{\flat}(\delta A)
\le M_{\flat}(A)
\mbox{ and }   E_{\flat}(\mu A)
\le 2M_{\flat}(A) $  for every flat chain $A$  (\cite{soapfilm} 5.2)).     

Henceforth we set $G = \Z_2$  and $n = 3$ for our application to soap films in three-space. In this case the weight of a $k$-dipolyhedron coincides with its Hausdorff $k$-measure.  Thus, for $k = 2$, the quantity $W(D)$ can be thought of as the {\em area} of $D$ and for $k = 1$ it is the {\em length} of $D$.  

  \section{Structure of dipolyhedra}
 
\subsection*{Splittings}  The orthogonal projection of a dipole $k$-chain $S$ into $\R^3 $      is a  mod two $k$-polyhedron  denoted  $\overline{S}$.  Then $S = \d
\overline{S}.
$  Similarly, the orthogonal projection of a mass $k$-chain $T$ into $\R^3$ is a  mod two $(k-1)$-polyhedron  denoted  $\overline{T}$ and satisfies 
$T =
\mu
\overline{T} .$  

\begin{prop}\label{converge} Suppose $ D = S+T \in \mathbf{D}_k(\Z_2)$.  Then  $M_{\flat}(\overline{S}) \le E_{\flat}(D); M_{\flat}(\overline{T}) \le E_{\flat}(D).$
\end{prop}

\begin{proof}  Given $\e > 0$ there exist $Q \in {\mathbf D}_k(\Z_2), R \in {\mathbf D}_{k+1}(\Z_2)$ with $D = Q + \p R$ and $E_{\flat}(D) > E(Q) + E(R) - \e = W(S_Q) + W(S_R) + M(T_Q) + M(T_R) -\e.$    Note that 
 $S_D = S_Q + S_{\p R}.$ Since $R = S_R + T_R$ then $\p R = \p S_R + \p T_R.$   Since the boundary of a dipole chain is also a dipole chain 
$S_{\p R} = \p S_R + S_{\p
T_R}$ 
 and therefore
$\overline{S} = \overline{S_D} =\overline{S_Q} + \overline{\p S_R} + \overline{S_{\p T_R}} = \overline{S_Q} +  \p \overline{  S_R} + \overline{T_R}.$
It follows that $$M_{\flat}(\overline{S}) \le  M(\overline{S_Q}) + M( \overline{  S_R}) + M( \overline{T_R}) = W(S_Q) + W(S_R) + M(T_R) < E_{\flat}(D) +
\e.$$

For the second inequality, note that $\p S_R$ is a dipole chain and 
$\p T_R$ is a sum of a dipole chain and a mass chain  $  T_{\p R}.$  As before, $\p R =  \p S_R + \p T_R.$   Therefore $\overline{T_{\p R}} = \p \overline{T_R}.$  
Since $\overline{T} = \overline{T_D} =\overline{T_Q} + \overline{T_{\p R}}$    
it follows that 
 
 $$M_{\flat}(\overline{T}) \le M(\overline{T_Q}) + M(\overline{T_R}) = M(T_Q) + M(T_R) < E_{\flat}(D) + \e.$$ 
   
   Since these inequalities hold for all $\e > 0$, the proposition follows.

  \end{proof}

\begin{thm}\label{split} If $A$ is a flat $k$-dipolyhedron then there exist a unique flat $k$-chain $B$ and a unique flat $(k-1)$-chain $C$ such that $A = \d B + \mu C$ and 
$E(A) = M(B) + M(C).$   If $D_i = \d B_i + \mu C_i \buildrel E_{\flat} \over \to A$ with $E(D_i) \to E(A)$ then $B_i \buildrel M_{\flat} \over \to  B, C_i \buildrel M_{\flat} \over
\to C$, $M(B_i) \to M(B)$ and
$M(C_i)
\to M(C).$
\end{thm}

\begin{proof}  Suppose $D_i \buildrel E_{\flat} \over \to A $ where $D_i = S_i + T_i$.   By Proposition \ref{converge} $B_i = \overline{S_i}$ is a Cauchy sequence in
$M_{\flat}.$  Let
$B$ denote its flat chain limit.  By 5.2 of \cite{soapfilm} $ S_i  =\d B_i  \buildrel E_{\flat} \over \to \d B.$   Similarly $C_i = \overline{T_i}$ converges
to a flat chain
$C$ in
$M_{\flat}$ and
$T_i = \mu C_i \buildrel E_{\flat} \over \to \mu C.$  It follows that $A = \d B + \mu C.$  

We next prove uniqueness:  Suppose $A = \d B + \mu C = 0.$  Then $E_{\flat}(A) = 0.$   By Propotision \ref{converge} $$M_{\flat}(B) = \lim M_{\flat}(B_i) \le \lim E_{\flat}(D_i) =
E_{\flat}(A) $$ and $$M_{\flat}(C) = \lim M_{\flat}(C_i) \le  \lim E_{\flat}(D_i) = E_{\flat}(A).
$$
It follows that $B = C = 0$  since $M_{\flat}$ is a norm.  Uniqueness of $B$ and $C$ follows.

Suppose $E(D_i) \to E(A).$  By lower semicontinuity of mass in the flat norm \cite{whitney},  $E(A)= M(B) + M(C) \le \liminf M(B_i) + \liminf M(C_i) \le \liminf E(D_i) = E(A).$   Thus $E(A) = M(B) + M(C) =
\liminf M(B_i) +  \liminf M(C_i).$  Since
$M(B) \le \liminf M(B_i)$ and $M(C) \le \liminf M(C_i)$ and all terms are nonnegative the result follows.
 
\end{proof}  
 
If $A = \d B + \mu C$ we say that $A$ {\em splits} into $\d B$ and $\mu C$. Since the splitting is unique  we may define  $W(A) = M(B)$ and $M(A) = M(C).$

  \begin{lem} \label{boundary}  If $A = \d B + \mu C$ is a flat dipolyhedron then $\p A = \d (\p B + C) - \mu (\p C).  $
\end{lem}

\begin{proof}  Acccording to (\cite{soapfilm}, 3.2)  $\p \d = \d \p$ and  $\d = \p \mu + \mu \p.$  The result follows .
\end{proof}

\begin{cor} \label{start} If $\g$ is a flat $(k-1)$-chain and  $A = \d B + \mu C$ is a flat $k$-dipolyhedron satisfying $\p A = \d \g$ then   $\p C = 0$ and $\p B + C = \g.$
\end{cor}

\begin{proof}  It follows from Lemma  \ref{boundary} that $\p A = \d(\p B +  C ) - \mu \p C = \d \g.$  Since the splitting is unique (Theorem \ref{split}) it follows that $\p B + C = \g $ and $\p
C = 0.$
\end{proof}
 \subsection*{The support of a flat dipolyhedron}   
 Let $B$ be a flat chain with finite mass.  According to   
(\cite{Fleming2}, \S 4) there exists a Borel measure $\rho_B$ and for every Borel set $X \subset
\R^3$ there exists a flat chain $B \cap X$ such that $\rho_B(X) = M(B \cap X).$     Moreover, if $P_i \buildrel M_{\flat} \over \to B$ with $M(P_i) \to
M(B)$ then
$P_i \cap X \buildrel M_{\flat} \over \to B \cap X $ and $M(P_i \cap X ) \to M(B \cap X )$ for every $X$ such that $\rho_B(fr X) = 0.$  (\cite{Fleming2} \S 4)
Define $(\d B) \cap X =\d(B \cap X)$ and $(\mu B) \cap X = \mu(B \cap X)$ the {\em part} of $\d B$ in $X$ and the {\em  part } of $\mu B$ in $X$, respectively.

Suppose $A \in  \mathcal{D}_k $ is a flat dipolyhedron with  bounded energy $E(A) \le \l$.    From Theorem \ref{split} we know $A = \d B + \mu C$ where $B$ and $C$ are flat chains with
$M(B) +M(C)
\le \l$.     Define  $$A \cap X = \d(B \cap X) + \mu(C \cap X).$$  We   call $A \cap X$ the {\em part of} $A$ in $X$.  
Define Borel measures  $\o_A(X) = \rho_B(X)$, $\mu_A(X)  = \rho_C(X),$  
and $\nu_A(X) =
\o_A(X) +
\mu_A(X).$
 The next proposition follows directly.
 
\begin{prop}\label{exceptional}  If $A = \d B + \mu C$  is a flat dipolyhedron with   finite energy and $X$ is a Borel set  there exists
a unique  flat dipolyhedron $A \cap X$ such that $W(\d B \cap X) =
\o_A(X)$,  
$ W(\d B -  \d B
\cap X) = \o_A(X^c),$ $M(\mu C \cap X) =
\mu_A(X)$,  
$ M(\mu C - \mu C
\cap X) = \mu_A(X^c)$, 
   $E(A \cap X) =
\nu_A(X)$ and  
$ E(A - A
\cap X) = \nu_A(X^c).$  If $D_i \buildrel E_{\flat} \over \to A$ with $E(D_i) \to E(A)$ then $D_i \cap X \buildrel E_{\flat} \over  \to A \cap X$ and $E(D_i \cap X) \to E(A \cap X)
$ for all $X$ such that
$\nu_A(fr X) = 0.$  
     
\end{prop}

     The {\em support} of a Borel measure $\nu$ is the smallest closed set $X$
whose complement is $\nu$-null.    Say that $\nu$ is a {\em measure on $Y$} if $Y$
contains the support of
$\nu$.

  We say a closed set
 $F$ {\em supports} a flat dipolyhedron  
 $A$ if for every open set
 $U$ containing $F$ there is a
 sequence $\{D_i\}$ of  dipolyhedra tending to $A$ in $E_{\flat}$ such
 that   $| D_i|\subset U$ for each $j$.  If there is a smallest set $F$
 which supports $A$ then $F$ is called the {\em support} of $A$ and
 denoted $| A|$.  The next theorem shows that every flat dipolyhedron with finite energy has a well defined support.  
 
 \begin{thm}\label{support}  If $A \in \mathcal{D}_k$ with $E(A) < \i$ then $|A|=
 |\nu_A|.$
 \end{thm}
 
 \begin{proof} By Theorem \ref{split} we know $A = \d B + \mu C$ where $B$ and $C$ are flat chains with $M(B) < \i$ and $M(C) < \i.$    Apply (\cite{Fleming2}, 4.3) to deduce   $|A|
 = |B| \cup |C| = |\rho_B| \cup  |\rho_C| = |\nu_A|.$
  \end{proof}

\section{Cones, pushforwards and projections}
 
\subsection*{Cones over dipolyhedra}   Let $\s$ be a $k$-cell in supported in $Q(p,r),$ the $3$-cube in $\R^3$
centered at $p$ with side length $r$.  The cone
$p
\s$ is also a cell found by intersecting the cones over the half-spaces forming $\s.$  Its boundary satisfies
$\p p\s = \s - p\p \s $ and $M(p\s) \le \frac{r \sqrt{3}}{k+1}M(\s).$  Define $p \d \s = \d p \s$ and $p \mu \s = \mu p \s.$
 
 It follows that  $p\d  \s$ is a well defined dipole cell with $W(  p \d \s) \le \frac{r\sqrt{3}}{k+1}W(\d  \s),$ and
$p\mu \s$ is a well defined mass cell with $M( p \mu \s) \le \frac{r\sqrt{3}}{k+1} M(\mu  \s).$     Extend the definition by
linearity to define dipole chains
$pS$ and mass chains
$pT$ taken over dipole chains
$S$ and mass chains
$T$, respectively.  Next extend the definition to dipolyhedra $pD = pS + pT.$ Observe $W(pS) \le \frac{r\sqrt{3}}{k+1}W(S)$ and $M(pT) \le
\frac{r \sqrt{3}}{k+1}M(T)$.

\begin{prop} \label{coneboundary}   If  $D \in {\mathbf D}_k $ then $pD \in{\mathbf D}_k$ and $$D = \p (p
D)  + p (\p D).$$  If  $D$   is
  supported in $Q(p,r)$    then  so is $pD$ and $$E(pD)
\le   \frac{r\sqrt{3}}{k+1} E(D).$$
\end{prop}

\begin{proof}
 
The first part reduces to showing $\d \s = \p(p\d \s)  + p \p \d \s $ and $\mu \s= \p(p\mu  \s)   + p \p \mu \s .$
 These follow directly from the definitions and linear relations.
  Suppose $D = S + T.$  Then  $pD  = pS + pT$ where $pS$ is a dipole chain and $pT$ is a mass chain.  Thus
  $$E(pD) = W(pS) + M(pT) \le \frac{r\sqrt{3}}{k+1}(W(S)+ M(T))  =  \frac{r\sqrt{3}}{k+1}E(D).$$
\end{proof}

\begin{thm} \label{Est}  If  $D    \in {\mathbf D}_k $   is
  supported in $Q(p,r)$    then  $$E_{\flat}(pD)
\le  \left(1+  \frac{r\sqrt{3}}{k+1}\right)E_{\flat}(D).$$
\end{thm}

\begin{proof} Let $\e > 0.$   There exist dipolyhedra $Q$ and $ R$ such that $D = Q + \p R$ and  $E_{\flat}(D) > E(Q) + E(R) -
\e.$  

By  Proposition \ref{coneboundary}  $p \p R =    R + \p(pR)$  for all dipolyhedra $R$. Then
$$E_{\flat}(p
\p R)
\le  E(R) + 
E(pR) \le 
 \left(1+  \frac{r\sqrt{3}}{k+1}\right)E(R).$$
  
Since $pD = pQ + p \p R$ it follows that
$$
\begin{array} {rll} E_{\flat}(pD) &\le E(pQ) +
E_{\flat}(p \p R) 
\\&\le   \left(  1+\frac{r\sqrt{3}}{k+1}\right)(  E(Q) + E(R)  )\\&\le
 \left(1+
\frac{r\sqrt{3}}{k+1}\right) (E_{\flat}(D) + \e).
\end{array}
$$   The result follows since this holds for all $\e > 0$.  
\end{proof}

It follows  that if $A$ is a flat dipolyhedron then the cone $pA$ has
unique definition as a flat dipolyhedron as follows:  if $D_i \buildrel E_{\flat} \over \to A$ then $\{pD_i\}$ forms a Cauchy sequence.  Denote its limit by $pA.$    
\begin{prop} \label{coneboundaryA}   If  $A \in \mathcal{D}_k(\Z_2) $ then  $$A = \p (p
A)   + p (\p A).$$  If  $A$   is
  supported in $Q(p,r)$    then  $$E_{\flat}(pA)
\le  \left(1+  \frac{r\sqrt{3}}{k+1}\right)E_{\flat}(A)$$ and $$E(pA)
\le   \frac{r\sqrt{3}}{k+1} E(A).$$
\end{prop}

\begin{proof}   The first two relations follow from Proposition \ref{coneboundary}, Theorem \ref{Est} and continuity of the boundary operator.
 Since energy is lower semicontinuous, there exists
$D_i
\buildrel E_{\flat}
\over
\to A$ such that
$E(D_i)
\to E(A).$   By Theorem \ref{Est}
$pD_i \buildrel E_{\flat} \over \to pA$ and hence
$$E(pA) \le  \liminf E(pD_i) \le \liminf  \frac{r\sqrt{3}}{k+1} E(D_i) =  \frac{r\sqrt{3}}{k+1}E(A).$$ 
\end{proof} 

\subsection*{Lipschitz pushfoward}  Let $f:U \subset \R^3 \to \R^3$ be a Lipschitz mapping.   Extend $f$ to $\R^{4}$ by $f(x,t) = (f(x),t).$  If $B$ is a flat $k$-chain supported
in an open set
$U \subset \R^3$ then the pushforward
$f_*B$ is well defined as a flat $k$-chain and satisfies $M_{\flat}(f_*B) \le |f|_{Lip}^kM_{\flat}(B)$ and $M(f_*B) \le |f|_{Lip}^k M(B)$  (\cite{whitney}).    It is called a  {\em Lipschitz chain}.   A chainlet $A = \d B +
\mu C$ is called a {\em Lipschitz dipolyhedron} if $B$ and $C$ are Lipschitz chains. If $D = \d B+ \mu C$ define $f_*D = \d f_*B + \mu f_* C.$

  \begin{prop}\label{chain} If $D  = \d B + \mu C $ is a  $k$-dipolyhedron and $f:U \subset \R^3 \to \R^3$ is a Lipschitz mapping  with $|D| \subset U$ then  $f_* D$ is a
$k$-dipolyhedron with  $\p f_*D  = f_* \p D$, $M(f_*D) \le |f|_{Lip}^k M(D), W(f_*D) \le |f|_{Lip}^k W(D)$ and $E_{\flat}(f_*D) \le |f|_{Lip}^k E_{\flat}(D).$
 
\end{prop}

\begin{proof} By Lemma \ref{boundary} we know that $\p f_* D = \d(\p f_*B + f_*C) - \mu (\p f_* C).$  Since $\p D = \d(\p B + C) - \mu (\p D)$ we have $f_*(\p D) = \d(f_*\p  B + f_*C) - \mu(f_* \p D)$.  Since $f_*$ is a chain map on flat chains we conclude $f_*(\p D) = \p f_* D$.  

By  Theorem \ref{split} and Proposition  \ref{chain}  it follows that $E(f_*D) = M(f_*B) + M(f_*C) \le |f|_{Lip}^k(M(B)+M(C)) = |f|_{Lip}^k(E(D)).$

 Let $\e > 0$.  There exists $D = Q + \p R$ with $E_{\flat}(D) > E(Q) + E(R) -\e.$  Since $f_*D = f_*Q + f_* \p R  = f_*Q + \p f_* R$   it follows that 
$E_{\flat}(f_*D) \le E(f_*Q) + E(f_*R) \le  |f|_{Lip}^k(E(Q) + E(R)) \le |f|_{Lip}^k(E_{\flat}(D) + \e).$

Finally, $W(f_*D) = M(f_*S) \le |f|_{Lip}^kM(S) = |f|_{Lip}^kW(D).$
\end{proof} 

Let $A$ be a flat dipolyhedron. By lower semicontinuity of energy we may choose $D_i \to A$ such that $E(D_i) \to E(A).$  It follows that $f_*A$ is a well defined  flat
dipolyhedron with  $E(f_*A) \le |f|_{Lip}^k E(A),$ and $E_{\flat}(f_*A) \le |f|_{Lip}^k E_{\flat}(A),$
\subsection*{Projection into a cube}    

For $x = (x^1, x^2, x^3) \in \R^3$, define $\|x\| = \max\{|x^1|, |x^2| , |x^3|\}.$    
   For  $r > 0$,   define 
 
 $$ f_r(x) = \begin{cases} x, &\|x\| \le r,\\ rx/\|x\|, &\|x\| > r. \end{cases}
 $$

 Observe that $f$ has Lipschitz constant $\le 1.$

Denote $D(r)= f_{r*}D$, the projection of a dipolyhedron $D$.  Since $\p(D(r)) = (\p D)(r)$   we can write $\p D(r)$ with ambiguity.

  Let $Q_r = Q(0,r).$ Then $f_{r*}(\R^3) = Q_r.$   
It follows from Proposition  \ref{chain}   projections $A(r)$  are uniquely
defined for all flat dipolyhedra $A$ with
 $E_{\flat}(A(r)) \le E_{\flat}(A) $  and  $E(A(r)) \le E(A).$

 Define $$\mathcal{B}_k(\Z_2)    = \{A \in \mathcal{D}_k(\Z_2) : E(A) + E(\p A) <
\i  \} \mbox{ and}$$    $$\mathcal{B}_k^0(\Z_2)    = \{A \in \mathcal{D}_k(\Z_2) : E(A) + E(\p A) <
\i,   |A| \mbox{ is compact } 
  \}.$$

 \section{A deformation theorem for flat dipolyhedra}

 The next result is the deformation theorem, first proved  for integer coefficients by Federer and Fleming \cite{FF}.   It  was extended to $Z_2$ coefficients in   \cite{ziemer} and to
abelian groups  in (\cite{Fleming2}, 7.3).

Let $\chi$ be an $\e$-cubical grid of $\R^3.$   A $k$-polyhedron $P$ is a  {\em polyhedron  of $\chi$} if $P$ is supported in the $k$-skeleton of $\chi$ and  $\p P$ is supported
in its $(k-1)$-skeleton.  
A $k$-dipolyhedron $D = \d B + \mu C$ is a dipolyhedron {\em of $\chi$} if $B$ and $C$ are both polyhedra of $\chi$.
 
\begin{thm}  \label{deformfleming} There exists a positive number $c = c(k,n)$ with the following property.  Given  $A \in
\mathcal{N}_k^0(\Z_2)$  and
$\e > 0$ there exist an $\e$-cubical grid $\chi$,  a polyhedral
$k$-chain $P$ of $\chi$,
$Q \in
\mathcal{N}_k^0(\Z_2) $ and 
$ R \in
\mathcal{N}_{k+1}^0(\Z_2)$ such that:
\begin{enumerate}
\item $A - P =  Q + \p R;$
\item $M(P) \le c(M(A) + \e M(\p A)),\\
M(\p P) \le cM(\p A), \\ M(Q) \le c\e M(\p A), \\M(R) \le c \e M(A);$
\item $|P| \cup |R| \subset 6\e$-neighborhood of $|A|,$\\
$|\p P| \cup |Q| \subset 6\e$-neighborhood of $|\p A|.$
\item If $A$ is a polyhedral (Lipschitz) chain then $Q$ and $R$ are polyhedral (Lipschitz chains).    
\end{enumerate}
\end{thm}

The proof makes use of successive radial projections from  well chosen points in the interior of each $k$-cube of $\chi$ onto its lower dimensional skeleton which minimize  distortion.   (See \cite{FF}.)  We refer to the {\em projection path} of a flat chain $A$
 as the rays traced by the projections of $|A|$. 
 
\begin{cor} \label{deformL} Let $\g$ be a Lipschitz Jordan curve in $\R^3$.  There exists a positive number $c = c(k,n)$ with the following property.  Given  $A \in
\mathcal{B}_k^0(\Z_2)$  with $\p A = \d \g$ and
$\e > 0$ there exist an $\e$-cubical grid $\chi$,  a dipolyhedral
$k$-chain $D$ of $\chi$,
$Q \in
\mathcal{B}_k^0(\Z_2)$ and  
 $R \in
\mathcal{B}_{k+1}^0(\Z_2)$
such that   
\begin{enumerate}
   \item $A - D = Q + \p R$;
\item $E(D)  \le 2c(E(A) + \e M(\p B))$,\\
$E(\p D) \le c M(\g), E(Q) \le c \e M(\g)$,\\
$E(R) \le c \e E(A)$;
\item $|D| \cup | R| \subset 6\e-$neighborhood of $| A|$,\\
  $| \p D| \cup |Q| \subset 6\e-$neighborhood of $| \p A|$.
\item If $A$ is a (Lipschitz) dipolyhedron then $Q$ and $R$ are (Lipschitz) dipolyhedra.  
 
\end{enumerate}

\end{cor}

\begin{proof} 
  
Suppose $A = \d B + \mu C$ with  $\p A = \d \g.$  By Corollary \ref{start} $\p C = 0$ and $\p B + C = \g$.  Apply Theorem \ref{deformfleming} to $B$ and $C$ to find $B- P_B = Q_B + \p R_B$ and $C - P_C = Q_C + \p
R_C$ satisfying properties (1)-(4) of Theorem \ref{deformfleming}. Therefore  $Q_C = 0$ and $\p P_C = 0.$   
Since $\g = \p B + C$ it follows that 
 
$$\g - (\p P_B  + P_C) = \p(Q_B + R_C).$$  Now the polyhedra $P, Q$ and $R$ are found by projecting $B$ and $C$ onto the $k$-skeleton of $\chi$.   
Thus $Q_B$ is the projection path of $\p B$, $R_C$ is the projection path of $C$.  Hence $Q_B + R_C$ is the projection path of $\p B + C$ which is the same as the projection path of $\g$.  Hence $$M(Q_B + R_C) \le c \e M(\g).$$ Similarly, $P_C $ is the projection of $C$, $\p P_B $ is the projection of $\p B;$ hence $P_C + \p P_B$ is the projection of $C + \p B$.  Hence   $$M(\p P_B + P_C) \le c M(\g).$$

Let $D = \d P_B + \mu P_C$.  Then $A - D = Q + \p R$ where $Q = \d(Q_B + R_C)$ and $R = \d R_B - \mu R_C.$  This establishes (1).  

By Theorems  \ref{split} and \ref{deformfleming} (2)
$$
\begin{array}{rll} E(D) = M(P_B) +
M(P_C) &\le c(M(B) + M(C) + \e (M(\p B) + M(\p C))) \\&\le c(E(A) + \e M(\p B)).
\end{array}
$$
Now $\p D = \d(\p P_B + P_C)$ since $\p P_C = 0.$  Then 
$$
\begin{array}{rll} E(\p D) = M(\p P_B + P_C) \le c M(\g).
\end{array}
$$
and
$$
\begin{array}{rll}  E(Q) = M(Q_B + R_C) \le c \e M(\g).
\end{array}
$$

Finally, $$E(R) = M(R_B) + M(R_C) \le c \e(M(B) + M(C)) = c \e E(A).$$   This completes the proof of (2).

The first part of (3) follows from the flat chain analogue since $|D|   = |P_B| \cup |P_C|$ and $|R| = |R_B| \cup |R_C|.$   For the second part, $Q$ is the projection path of $\p
A$.  Thus $|\p Q| \subset 6 \e-$nbd of $|\p A|.$   

Part  (4) is an easy consequence of the flat chain analogue and the definitions of $A, Q$ and $R$.

\end{proof}

  We say that a flat dipolyhedron $A$  {\em spans} $\d \g$ if  $\p A = \d \g$ and  the following condition holds: if $X$  is a $2$-dimensional subspace of $\R^3$  and  $\Pi: \R^3 \to X$ is an orthogonal  projection which is an immersion of $\g$, then  $\p \Pi_*A = \Pi_* \p A.$ 
 
 \begin{figure}[b]
\begin{center}
\resizebox{4.0in}{!}{\includegraphics*{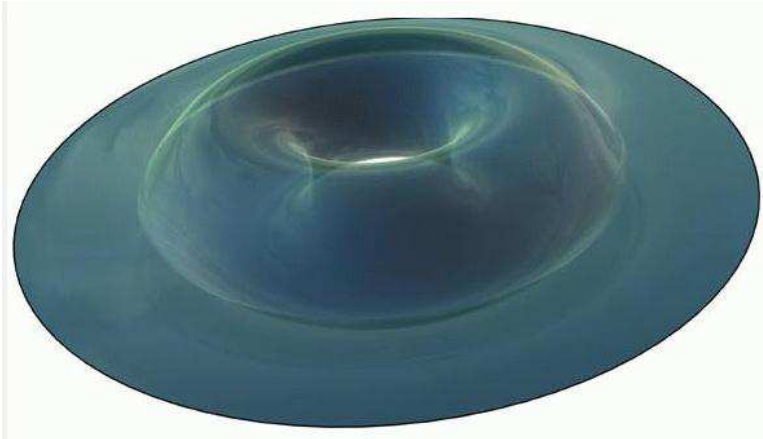}}\label{test}
\caption{A dipolyhedron that does not span its boundary. \tiny{ Drawing by Harrison Pugh}}
\end{center}
\end{figure}

\begin{lem}  If $D_i \to A$   and $D_i$ spans $\d \g$ then $A$ spans $\d \g$.  
\end{lem}

\begin{proof}  By continuity of the boundary and pushforward operators  $\p D_i \to \p A = \d \g$ and 
$\p \Pi_*A = \lim \p \Pi_* D_i = \lim \Pi_* \p D_i = \Pi_* \p A.$
\end{proof}

\begin{lem}\label{delta} There exists $\e > 0$ such that  if  $A$ spans $\d \g$ then $W(A) > \e.$
\end{lem}

\begin{proof}  Choose a projection $\Pi$ so that $\Pi_* \g$ is a Jordan curve in $X$.  Let $K$ denote the chain whose support is the inside of  $\Pi_* \g$ in $X$ so that $\p K = \Pi_* \g.$  Let $\e = M(K).$     The condition that  $\p \Pi_*A = \Pi_* \p A = \Pi_* \g$ implies that $ \Pi_*A = K$.  The result follows since $W(A) \ge W(\Pi_*A) \ge \e.$  \end{proof} 

 Choose $\l$ sufficiently large so that $\g \subset  Q_{\l'}$ where $\l' = \l(k+1)/M(\g).$ 

Define $$\G(\l) = \{ A \in \mathcal{B}_k^0(\Z_2):  E(A) \le \l, |A| \subset Q_{\l'}, A \mbox{ spans } \d \g\}.$$    This collection of supports $|A|$ of dipolyhedra in $\G(\l)$ contains
solutions to other Plateau type problems including the supports of
\begin{itemize}
\item  mappings of the two disk $B$,
$f:B \to \R^3$, $f(\p B) = \g$ and $f$ is smooth away from $\p B$.  
\item area minimizing integral currents
 
\item area  minimizing soap films $S$, as observed by Plateau:     The set $S$ is smooth away from its
branch set $B$.    Consider the connected components
$X_i$ complementary to $B$.  These are embedded and smooth mod two surfaces.  Let $X = \sum X_i$     and $D = \d X + \mu B .$
Then $|D| = |X| = S$  and $\p D = \d \g.$ 
\end{itemize}   
 
\begin{thm} \label{compact} $\G(\l)$ is compact and nonempty in the $E_{\flat}$ norm.  
\end{thm}

\begin{proof}  Suppose $A_i \buildrel E_{\flat} \over \to A$ where $A_i \in \G(\l).$  We know $E(A) \le \l$ by
lower semicontinuity of energy,  and $|A|
\subset Q_{\l'}$  since each $|A_i|
\subset Q_{\l'}.$  Since each $A_i$ spans $\d \g$ it follows that $A$ spans $\d \g$.  Thus $\G(\l)$ is closed.  Use Corollary \ref{deformL} to show it is totally bounded.  Given $\e > 0$ there exists
a dipolyhedron $D$ such that $A - D = Q + \p R$ with 
$$E_{\flat}(A-D) \le E(Q) + E(R) \le c \e(E(\p A) + E(A)) \le c \e (M(\g) + \l).$$  It follows that $\G(\l)$ is totally bounded and thus compact.  
 
We show that $\d 0\g \in
\G(\l).$  Since $0 \g$ is a cone over a polyhedron it is a flat chain with $E(\d 0 \g) = M(0 \g) \le \frac{\l'}{k+1} M(\g) = \l$    ( \cite{Fleming2} \S 6).
 By Proposition \ref{coneboundary} $\p \d 0\g = \d \g $ and we know $\d 0\g$   spans its
boundary.    Since $|\d 0\g| \subset Q_{\l'}$ the
  result follows
\end{proof}

 Let $m = \inf\{W(A): A \in  {\G}(\l)\}.$  There exist
$A_i
\in \G(\l)$ such that
$W(A_i)
\to m.$  By compactness (Theorem \ref{compact}) the sequence
$A_i$ has a subsequential limit
$A
\in \G(\l)
$ with
$W(A) \le \liminf{W(A_i)} \le m.$  Then $W(A) = m.$  According to Lemma \ref{delta} it follows that $m > 0.$

  It may be that there is another $A$ with smaller area spanning $\g$ with $E(A) \le\l$ that
 is not supported in
$Q_{\l}$.   Let $f_{\l}$ denote the projection into $Q_{\l}$.  Since $A$ spans its boundary so does $A_{\l} = f_{\l}A$.   From   Proposition \ref{chain} we
conclude 
$A_{\l} 
\in
{\G}(\l)$.  By Proposition  \ref{chain}  $W(A_{\l} ) \le W(A)$  and $E(A_{\l} ) \le E(A) \le \l$,  so we may replace $A$ with
$A_{\l} $.   
 
The flat dipolyhedron $A$ is our   solution to Plateau's problem for soap films with   energy bounded by $\l$.    

\begin{thm}[Almost everywhere regularity] The set $ |A|\backslash |\p A|$ is a smooth surface except on a union of Lipschitz Jordan curves with finite total length.
\end{thm} 

\begin{proof}  By Theorem \ref{split} $A = \d B + \mu C$ where $B$ and $C$ are flat chains (mod two) with finite mass and
finite boundary mass.      Decompose B into its indecomposable parts $B = \sum B_i$. ) Since $B_i$ is
area minimizing, we may apply mod two  regularity (\cite{Fleming1}) to deduce each surface
$|B_i|$   is smoothly embedded away from its boundary.  
  In his thesis, Ziemer proved that mod two boundaries are integral currents (\cite{ziemer},
6.5).   Fleming proved that
  integral 1-cycles with finite mass are  sums of closed
curves each of which is Lipschitz.   Since $\p C = 0$ and $M(C) < \l$ then $C$ is a sum of Lipschitz Jordan curves  with finite total length.   (See \cite{federer},  4.2.25)\end{proof}


\begin{thebibliography}{rll}
 \bibitem[A]{A} 	F.J. Almgren, Jr.,  Plateau's Problem,
Benjamin, New York, 1966.
\bibitem[D]{Douglas}		Jesse Douglas, Solution of the
problem of Plateau,
		Trans. Amer. Math. Soc. 33 (1931), 263-321.

 
 \bibitem[Fe] {federer} Federer, Herbert, {\sc Geometric Measure Theory},
Springer-Verlag, New York, 1969

\bibitem[FF]{FF}		Herbert Federer and Wendell H. Fleming,
		Normal and integral currents, Ann. of Math. 72
	(1960), 458-520.

\bibitem[Fl1]{Fleming2} Wendell H. Fleming, Flat chains over a finite
coefficient group, Trans. AMS,
Jan-Feb 1966 v 121 pp160-186

\bibitem[Fl2]{Fleming1}		Wendell H. Fleming, On the oriented Plateau
		problem, Rend. Circ. Mat. Palermo (II) 11 (1962),
 		1-22

\bibitem[H]{soapfilm} Harrison, Jenny,  Cartan's magic formula and soap film structures,
to appear, Journal of Geometric Analysis 

  
 
\bibitem[O]{Osserman} Robert Osserman,  A Proof of the Regularity Everywhere of the Classical Solution to Plateau's Problem,  Ann. of Math. 91, (1970), 550-569.
\bibitem[P]{Plateau}		J.A.F. Plateau, Statique Experimentale
		et Theorique des Diquides Soumis aux Seules
		Forces Moleculaires, Paris, Gauthier-Villars, 1873.

  
 \bibitem[W]{whitney} Hassler Whitney, Geometric Integration Theory,  Princeton University Press, 1957.

\bibitem[Z]{ziemer} William P. Ziemer,  Integral currents mod 2, Trans of the AMS, Vol 105, Issue 3,   (1962), 496-524

\end{thebibliography}
\end{document}